\documentclass[12pt]{article}

\usepackage[latin1]{inputenc}
\usepackage{amsmath,amssymb}
\usepackage{amsmath,accents}
\usepackage{alltt,graphics,graphpap,color}
\usepackage[dvips]{graphicx}
\usepackage{amsfonts}
\usepackage{amscd}
\usepackage{mathrsfs}
\usepackage{hyperref}
\usepackage{float}
\usepackage{subcaption} 

\DeclareMathAlphabet{\mathpzc}{OT1}{pzc}{m}{it}

\newtheorem{Theorem}{Theorem}[section]

\newtheorem{Lemma}[Theorem]{Lemma}
\newtheorem{Corollary}[Theorem]{Corollary}
\newtheorem{Remark}[Theorem]{Remark}

\newcommand{\mm}{{\cal M}}

\newcommand{\bg}{\bar{g}}
\newcommand{\proof}{\emph{Proof. }}
\newcommand{\cvd}{\hfill$\square$ \bigskip}

\newcommand{\sss}{Sol_3}
\newcommand{\cth}{\cos\vartheta}
\newcommand{\sth}{\sin\vartheta}

\newcommand{\sthq}{\sin^2\vartheta}
\newcommand{\vth}{\vartheta}
\newcommand{\vthz}{\vartheta_0}
\newcommand{\vthu}{\vartheta_1}
\newcommand{\vthd}{\vartheta_2}

\begin{document}

\def\qed{\hbox{\hskip 6pt\vrule width6pt height7pt
depth1pt  \hskip1pt}\bigskip}


\title{Invariant translators of the Solvable group}

\author{\sc Giuseppe Pipoli}
\date{}

\maketitle

{\small \noindent {\bf Abstract:} We classify the translators to the mean curvature flow in the three-dimensional solvable group $\sss$ that are invariant under the action of a one-parameter group of isometries of the ambient space. In particular we show that $\sss$ admits graphical translators defined on a half-plane, in contrast with a rigidity result of Shahriyari \cite{Sh} for translators in the Euclidean space. Moreover we exhibit some non-existence results.} 

\medskip


\bigskip

\section{Introduction}

The aim of this paper is to classify all the translators of the three-dimensional solvable group $\sss$ that are invariant under the action of a one-parameter group of symmetries. 

A hypersurface $\mm$ in a general ambient space is said a \emph{soliton} to the mean curvature flow if it evolves without changing its shape. More precisely it means that there exists $G=\left\{\left.\varphi_t\ \right|\ t\in\mathbb R\right\}$ a one-parameter group of isometries of the ambient space such that the evolution by mean curvature flow of $\mm$ is given by the family of hypersurfaces
\begin{equation}\label{GGGX}
\mm_t=\varphi_t\left(\mm\right).
\end{equation}
When the ambient space is a Lie group equipped with a left invariant metric, then we can choose $G$ consisting of left translations. In this case we say that $\mm$ is a \emph{translator}.  

Translators provide interesting examples of explicit eternal solutions of the flow and, as shown by Huisken and Sinestrari \cite{HS}, they have a fundamental role in the analysis of type II singularities.

It is well known that the property of being a translator can be view as a prescribed mean curvature problem: let $V$ be the Killing vector field associated to $G$, then $\mm$ evolves moving according to \eqref{GGGX} if and only if
\begin{equation}\label{defi}
H=\bar g(\nu,V),
\end{equation}
where $\nu$ and $H$ are respectively the unit vector field and the mean curvature of $\mm$, while $\bar g$ is the metric in the ambient space. In this case we will say that $\mm$ is a translator in the direction of $V$. A proof of \eqref{defi} can be found in \cite{HuSm}.

The ambient space considered in this paper is the three-dimensional solvable group $\sss$. It is one of the eight geometries of Thurston, the one with the least number of isometries. It can be seen as a Riemannian Lie group defined in the following way: on $\mathbb R^3$ equipped with the usual coordinates we consider the group operation,
$$
(x_1,y_1,z_1)\star(x_2,y_2,z_2)=\left(x_1+e^{-z_1}x_2,y_1+e^{z_1}y_2,z_1+z_2\right),
$$
and the left invariant Riemannian metric
$$
\bar g=e^{2z}dx^2+e^{-2z}dy^2+dz^2.
$$
It follows that an orthonormal basis of left invariant vector fields is
$$
E_1=e^{-z}\frac{\partial}{\partial x},\quad E_2=e^{z}\frac{\partial}{\partial y},\quad E_3=\frac{\partial}{\partial z}.
$$
The geometry of this manifold is well known, in particular its isometry group is generated by the left translation, hence it is $3$-dimensional. Moreover any Killing vector field is a linear combination with constant coefficients of 
$$
F_1=\frac{\partial}{\partial x},\quad F_2=\frac{\partial}{\partial y}, \quad F_3=-x\frac{\partial}{\partial x}+y\frac{\partial}{\partial y}+\frac{\partial}{\partial z}.
$$

Every time we wish to study hypersurfaces satisfying some curvature condition (like minimal, CMC, or constant Gaussian curvature), a successful strategy for finding examples is to look for the symmetrical ones. This idea was considered in $\sss$ \cite{Lo, LM, LM2}, but also in many other ambient spaces: without completeness claims, we mention \cite{FMP} for CMC in the Heisenberg group, \cite{On, MO1} in $\mathbb H^2\times\mathbb R$ and \cite{MO2} for more general ambient spaces. 

Since the mean curvature flow preserves the symmetries (see \cite{Pi1} for a proof), the same strategy can be applied for equation like \eqref{defi}. In fact the first known examples of translators in the Euclidean space are the tilted grim reapers cylinders (invariant by translation), the bowl solution and the translating catenoids \cite{CSS} (rotationally symmetric) and the translating helicoids \cite{Ha}. The author of the present paper described in \cite{Pi} the analogous of these examples in the Heisenberg group. It is not difficult to imagine that many other non-symmetric examples exist. Without completeness claim we mention some recent literature about the construction of interesting families of translators in the Euclidean space: \cite{BLT, HIMW, HMW}.

The goal of the present paper is to classify the translators in $\sss$ that are invariant under the action of a one-parameter group of isometries. Let $X$ and $V$ be two non-vanishing Killing vector field. With an abuse of notation we denote with $X$ also the one-parameter group of isometries associated to the vector field $X$. We want to describe the $X$-invariant translators in the direction of $V$. Since our ambient manifold is homogeneous we can always suppose that $(0,0,0)$ belongs to the surface considered. On the hoter hand, $\sss$ is not isotropic, then the result depends strongly on the interaction between $X$ and $V$.

The most interesting examples appear in the very special case $X=F_1$. We point out that the minimal surfaces with this kind of symmetry are classified by Lopez and Munteanu in \cite{LM}: they are the planes
\begin{equation}
y=y_0,\quad z=z_0,\label{piani}
\end{equation}
and the logarithmic surfaces
\begin{equation}
z=\log(y_0\pm y)+z_0,\label{log}
\end{equation}
where $y_0,z_0\in\mathbb R$. Clearly, by \eqref{defi}, these surface can be seen as trivial examples of translators when $V$ is tangent to them. For later use we introduce a further family of $F_1$-invariant surfaces that we can call of \emph{half-logarithmic} type:
\begin{equation}\label{half}
z=\frac 12 \log(y_0\pm y)+z_0,
\end{equation}
where once again $y_0,z_0\in\mathbb R$.

The whole family of $F_1$-invariant translators are described by our first result.

\begin{Theorem}\label{main1}
Let  $V=\eta F_1+\lambda F_2+\mu F_3$ be a non-zero Killing vector field and let $\mm$ be a $F_1$-invariant translator in the direction of $V$. Then we have one of the following possibilities.
\begin{itemize}
\item[1)] $\mm$ is one of the minimal surfaces \eqref{piani}, \eqref{log} if and only if $V$ is tangent to it.
\item[2)] If $\mu=0$ and $\lambda\neq 0$, then $\mm$ is a complete graph $y=F(x,z)$, for some function $F$ defined on a slab or on an half-plane. The lower half surface $\mm^-=\mm\cap\left\{z<0\right\}$ is always asymptotic to a horizontal plane, while the upper half surface $\mm^+=\mm\cap\left\{z>0\right\}$ is asymptotic to one of the following: a horizontal plane, a surface of logarithmic type \eqref{log}, a surface of half-logarithmic type \eqref{half} or a generic $F_1$-invariant plane  $y=C_1z+C_0$, for some constants $C_0,C_1$.
\item[3)] if $\mu\neq 0$, there is a line $l=\left\{(s,y_0,z_0)|s\in\mathbb R\right\}$ such that $\mm\backslash l$ has two connected components. Each one of them is asymptotic to the special vertical plane $y=-\frac{\lambda}{\mu}$ or to a minimal surface of logarithmic type \eqref{log}.
\end{itemize}
\end{Theorem}
This Theorem shows that in $\sss$ there are some graphical translators that are defined on a half-plane. This property cannot be satisfied for translators in the Euclidean space as proved in \cite{Sh}. Moreover some of them are not convex in contrast with another rigidity property of translators in the Euclidean space \cite{SX}. See Remark \ref{rmk01} for a more precise statement. This phenomenon has been observed in the Heisenberg group too \cite{Pi}. We point out that for almost all the cases considered in part $2)$ of Theorem \ref{main1} $\mm^+$ is asymptotic to a minimal surface. The latter two possibilities appear only in very special cases. We refer to Lemma \ref{f2f2} for more details. Some numerical simulations of the translators described in Theorem \ref{main1} can be seen in figures from \ref{fig01} to \ref{fig03} for part $2)$ and figures from \ref{fig04} to \ref{fig07} for part $3)$.

\noindent For a more general kind of invariance we do not have a precise description but we still have some structural properties.

\begin{Theorem}\label{main2}
Let $X=aF_1+bF_2$, let $V$ any Killing vector field and let $\mm$ be a $X$-invariant translator in the direction of $V$. If $ab\neq0$, then we have two possibilities:
\begin{itemize}
\item[1)] either $\mm$ is the horizontal plane $z=0$ (hence it is minimal) and $V$ is tangent to it,
\item[2)] or $\mm$ is a complete graph $y=F(x,z)$ for some function $F$.
\end{itemize}
\end{Theorem}

\noindent We complete the classification with the most rigid case.

\begin{Theorem}\label{main3}
Let $X=aF_1+bF_2+cF_3$, let $V$ be any Killing vector field and let $\mm$ be a $X$-invariant surface. If $c\neq 0$ the following are equivalent:
\begin{itemize}
\item[1)] $\mm$ is a translator in the direction of $V$, 
\item[2)] $\mm$ is minimal and $V$ is tangent to $\mm$.
\end{itemize}
\end{Theorem}

As a consequence of this Theorem we have that for suitable combinations of $X$ and $V$, there are no $X$-invariant translator in the direction of $V$. For example there are no such surfaces if $a=b=0$, $c=1$ and $\eta\lambda\neq 0$. See Remark \ref{rmkfin} for the proof. We mention that the $F_3$-invariant minimal surfaces are studied by Lopez in \cite{Lo}. 

Finally we point out that Theorems \ref{main1}, \ref{main2} and \ref{main3} exhaust all the possibilities, in fact the case $a=c=0$ is left, but it can be easily traced back to Theorem \ref{main1} applying to the whole ambient space the isometry $\phi:(x,y,z)\mapsto(x,y,-z)$.

The paper is organized as follows. In Section \ref{Prem} we collect some basic facts about the geometry of $\sss$, in particular we describe all the one-parameter group of isometries of this space. In Section \ref{SEZA} we start proving Theorem \ref{main1} and we exhibit some numerical simulations of these surfaces. We finish the paper with the proof of Theorem \ref{main2} and \ref{main3} in Section \ref{SEZB}.

\section{Preliminaries}\label{Prem}
The geometry of $\sss$ is well known, here we write only what is strictly necessary for this paper. Other information about this space can be found for example in \cite{Tr}. We racall that a left-invariant orthonormal basis is $$
E_1=e^{-z}\frac{\partial}{\partial x},\quad E_2=e^{z}\frac{\partial}{\partial y},\quad E_3=\frac{\partial}{\partial z}.
$$
The Levi-Civita connection of $\bar g$ is determined by
\begin{equation}\label{LC}
\begin{array}{lll}
\bar\nabla_{E_1} E_1=-E_3& \bar\nabla_{E_2}E_1=0 & \bar\nabla_{E_3}E_1=0\\
\bar\nabla_{E_1}E_2=0 & \bar\nabla_{E_2}E_2=E_3 & \bar\nabla_{E_3}E_2=0\\
\bar\nabla_{E_1}E_3=E_1 & \bar\nabla_{E_2}E_3=-E_2 &\bar\nabla_{E_3}E_3=0
\end{array}
\end{equation}

For every $p\in\sss$, we denote by $L_p:q\in\sss\mapsto L_p(q)= p\star q$ the left translation by $p$. The following result probably was already known, but we give the proof for completeness.

\begin{Lemma}\label{sottogruppi}
Let $X=aF_1+bF_2+cF_3$ be a non vanishing Killing vector field, then the one-parameter group of isometries of $\sss$ generated by $X$ is
$$
\begin{array}{ll}
\left\{\left.L_{\left(\frac ac(1-e^{-ct}),\frac bc (e^{ct}-1), ct\right)}\right| t\in\mathbb R\right\}, &\text{if}\ c\neq 0;\\
&\\
\Bigl\{\left.L_{\left(at,bt,0\right)}\right| t\in\mathbb R\Bigr\}, &\text{if}\ c=0.
\end{array}
$$
\end{Lemma}
\proof
Let $A,B,C$ be three real functions such that $\left\{\left.L_{\left(A(t),B(t),C(t)\right)}\right|t\in\mathbb R\right\}$ is a one-parameter group of isometries of $\sss$. Then $A(0)=B(0)=C(0)=0$ and for every $t,s\in\mathbb R$ we have
\begin{equation}\label{gruppo1}
\begin{array}{l}
\left(A(t),B(t),C(t)\right)\star\left(A(s),B(s),C(s)\right)\\
\phantom{aaaa} = \left(A(t)+e^{-C(t)}A(s),B(t)+e^{C(t)}B(s),C(t)+C(s)\right)\\
\phantom{aaaa} =\left(A(t+s),B(t+s),C(t+s)\right).
\end{array}
\end{equation}
Looking at the third components of \eqref{gruppo1}, we deduce that $C$ is a linear function, hence there is a constant $c$ such that $C(t)=ct$. If $c=0$, then we can see that $A$ and $B$ are linear too. Otherwise we can compute the derivative with respect to $s$ in $s=0$ of each member of the equality \eqref{gruppo1} and we get that for every $t$
\begin{eqnarray*}
A'(t) & = & A'(0)e^{-ct},\\
B'(t) & = & B'(0)e^{ct}.
\end{eqnarray*}
The two ODEs can be easily solved. Finally, deriving with respect to $t$, we can see that $X$ is the Killing vector field associated to the group.
 \cvd

\section{$F_1$-Invariant translators}\label{SEZA}

The goal of this section is to prove Theorem \ref{main1} describing the $F_1$-invariant translators. We recall that, as noticed in \cite{LM}, it is equivalent to study the $F_2$-invariant translators: in fact it is sufficient to apply to the whole ambient space the isometry $\phi(x,y,z)=(x,y,-z)$ to transform a surface of the first kind to one with the other direction of symmetry.

Let $\mm$ be a $F_1$-invariant surface, then there is a planar curve $\gamma(s)=(0,y(s),z(s))$ parametrized by arc length such that the surface is
\begin{equation}\label{rrr}
\mm(u,s)=(u,0,0)\star\gamma(s)=\left(u,y(s),z(s)\right).
\end{equation}
Following the idea of \cite{LM}, we define $\vartheta(s)$ as the angle between $\gamma'$ and $E_2$, hence we have the following ODE system.
\begin{equation}\label{sistema01}
\left\{\begin{array}{rcl}
y'&=&e^z\cos\vartheta;\\
z'&=&\sth;\\
\vartheta' & = & H;
\end{array}\right.
\end{equation}
where the derivatives are rispect to $s$ and $H$ is the mean curvature of $\mm$. Note that, because $\mm$ is invariant, $H$ depends only on $s$. Moreover, unlike \cite{LM}, we prefer to define the mean curvature as the sum of the principal curvatures (and not as the mean of the pricipal curvatures).
Let $V=\eta F_1+\lambda F_2+\mu F_3$ be a Killing vector field. Note that $F_1$ is tangent to $\mm$, therefore by the definition of translators \eqref{defi} we can consider $\eta=0$ without loss of generality. 
As shown in \cite{LM}, $\nu=-\sth E_2+\cth E_3$, while the Killing vector field can be rewritten as
$$
V=\lambda e^{-z}E_2+\mu\left(-xe^zE_1+ye^{-z}E_2+E_3\right).
$$
Therefore, by \eqref{defi}, $\mm$ is a translator in the direction of $V$ if and only if
\begin{equation}\label{H1}
H=-\lambda e^{-z}\sth+\mu\left(\cth-y e^{-z}\sth\right).
\end{equation}
By \eqref{sistema01} we get
\begin{eqnarray}
 \vth ' & = & -\lambda e^{-z}\sth+\mu\left(\cth-y e^{-z}\sth\right)\label{vth2}\\
\nonumber & = & -\lambda e^{-z}z'+\mu\left(e^{-z}y'-ye^{-z}z'\right)\\
& = & \frac{d}{ds}\left[e^{-z}\left(\lambda+\mu y\right)\right].\label{vth1}
\end{eqnarray}
Obviously the solution of \eqref{sistema01} is uniquely determined fixing the initial conditions. Since $\sss$ is homogeneous, we can always suppose without loss of generality that $y(0)=z(0)=0$, then for every fixed $\lambda$ and $\mu$, we have just one degree of freedom given by $\vth_0:=\vth(0)$. Once the initial conditions have been chosen, by \eqref{vth1} we get
\begin{equation}\label{vth3}
\vth-\vthz+\lambda=e^{-z}(\lambda+\mu y).
\end{equation}
Putting \eqref{vth3} into \eqref{vth2} we get that the third equation of the system \eqref{sistema01} becomes an autonomous independent equation:
\begin{equation}\label{vth4}
\left\{\begin{array}{rcl}
\vth ' & = & f(\vth):=\mu\cth-(\vth-\vthz+\lambda)\sth,\\
\vth(0) & = & \vthz.
\end{array}\right.
\end{equation}

The crucial part of the present paper is a careful estimation of the solution of the ODE \eqref{vth4}. After that the behaviour of $y$ and $z$ can be easily understood by \eqref{sistema01} and \eqref{vth3}. By periodicity it would be sufficient to consider $\vthz\in\left[0,2\pi\right]$, but we prefer to let $\vthz$ vary in the whole $\mathbb R$. Particular attention will be given to the zeros of the function $f$.

\begin{Lemma}\label{zero f}
For any $\lambda, \mu, \vthz\in\mathbb R$ the set of the zeros of $f$ is unbounded from above and from below. Moreover either $f(\vthz)=0$ or there are $\vthu,\vthd\in\mathbb R$ such that $\vthz\in\left(\vthu,\vthd\right)$, $f(\vthu)=f(\vthd)=0$ and $f$ has a sign in $\left(\vthu,\vthd\right)$. In particular we have that $\vthd-\vthu\leq\pi$ if $\mu=0$ and $\vthd-\vthu<2\pi$ otherwise.
\end{Lemma}
\proof When $\mu=0$ the equation $f=0$ is very easy to solve and the thesis is trivial. Let us consider the case $\mu\neq 0$. We have that for any $k\in\mathbb Z$
$$
f(k\pi)f((k+1)\pi)=-\mu^2<0.
$$
Hence the thesis follows because $f$ is a continuous function.\cvd

As an immediate consequence we have the following result.

\begin{Corollary}\label{vth monotona}
The function $\vth$, maximal solution of \eqref{vth4}, is defined for any $s\in\mathbb R$, is monotone and bounded. More precisely it evolves according one of the following mutually exclusive possibilities:
\begin{itemize}
\item[1)] if $f(\vthz)=0$ then $\vth(s)=\vthz$ for every $s$;
\item[2)] if $f(\vthz)>0$ then $\vth$ is strictly increasing and there exist the limits 
$$\lim_{s\rightarrow +\infty}\vth(s)=\vthd,\quad \lim_{s\rightarrow -\infty}\vth(s)=\vthu;$$
\item[3)] if $f(\vthz)<0$ then $\vth$ is strictly decreasing and there exist the limits 
$$\lim_{s\rightarrow +\infty}\vth(s)=\vthu,\quad \lim_{s\rightarrow -\infty}\vth(s)=\vthd.$$
\end{itemize}
\end{Corollary}

Since $\vth$ is defined for any $s$, by \eqref{sistema01} also $y$ and $z$ are defined on the whole real line too. Therefore we get that the curve $\gamma$ and the associated surface $\mm$ are complete. Moreover if $f(\vthz)=0$, then $H=\vth'=0$ and $\mm$ is a $F_1$-invariant minimal surface. In this case equation \eqref{defi} is satisfied if and only if $V$ is tangent to $\mm$. Such kind of minimal surfaces are classified in \cite{LM}: they are planes \eqref{piani} or logarithmic surfaces \eqref{log}. From now on we will consider always the case $f(\vthz)\neq 0$. As a consequence of these first results we have some restrictions on the behavior of the coordinate functions $y$ and $z$. 

\begin{Lemma}\label{punti critici}
If $\mu=0$, the function $z$ is strictly monotone, while $y$ has at most one critical point. If $\mu\neq 0$, the function $z$ has at most one critical point, while $y$ has at most two critical points. It follows that for any value of $\mu$, both the functions have a limit (finite or not) as $s$ goes to $\pm\infty$.
\end{Lemma}
\proof We proved in Corollary \ref{vth monotona} that $\vth$ is strictly monotone for any value of the parameters $\lambda$, $\mu$ and $\vthz$. Let us consider firs the case $\mu=0$: Lemma \ref{zero f} says that $\vth$ varies  in an interval of length smaller than $\pi$, then, by \eqref{sistema01}, we get that each one of the two functions can have at most one critical point. Moreover, since $\mu=0$, if $z'=0$ then $f=0$ too. Therefore $z$ cannot have a critical point in the open interval $(\vthu,\vthd)$ by definition of $\vthu$ and $\vthd$. When $\mu\neq 0$, Lemma \ref{zero f} says that $\vth$ varies  in an interval of length smaller than $2\pi$, therefore the critical points of the two functions are at most two for each. Finally, as shown in the proof of Lemma \ref{zero f}, at least one zero of $f$ lies between $k\pi$ and $(k+1)\pi$ for any $k\in\mathbb Z$. Therefore $z'=\sth=0$ can have at most one solution in the open interval $(\vthu,\vthd)$.\cvd

\subsection{$F_1$-invariant translators in the direction of $F_2$}
In this section we describe the properties of $\gamma(s)$ when $\mu=0$ and $\lambda\neq 0$. As a consequence we will prove part $2)$ of Theorem \ref{main1}. At first we can notice that, by Lemma \ref{punti critici}, $z$ is strictly monotone, hence the curve $\gamma$ and the surface $\mm$ are embedded and $\mm$ is a graph $y=F(x,z)$ for some function $F$. The following results describe the asymptotic behavior of $\gamma$.

\begin{Lemma}\label{main F2}
If $\mu=0$, for every $\vthz$ and $\lambda\neq 0$ we have that, as $|s|$ diverges,  the function $y(s)$ diverges while we have two possibilities for $z(s)$:
\begin{itemize}
\item[1)] $z(s)$ diverges to $+\infty$ if and only if $\lim_{|s|\rightarrow +\infty}\vth(s)=\vthz-\lambda$;
\item[2)] $z(s)$ converges to a finite constant otherwise.
\end{itemize}
In particular $z$ in bounded from below.
\end{Lemma}
\proof
For simplicity we consider only the case of $s\rightarrow +\infty$, the other one is analogous. Since $\mu=0$, \eqref{vth3} assume the simpler form
\begin{equation}\label{eee}
\vth-\vthz+\lambda=\lambda e^{-z}.
\end{equation}
By Corollary \ref{vth monotona} we know that the left hand term of \eqref{eee} converges to a finite constant. It follows that $z$ cannot be arbitrary negative. Furthermore $z$ converges to a finite constant if and only if the left hand term does not converge to zero, that is if $\lim_{s\rightarrow +\infty}\vth(s)\neq\vthz-\lambda$. 

Once we know the behavior of $z$, we can describe the evolution of $y$ too. Recalling that $y'=e^z\cth$, if $\lim_{s\rightarrow +\infty}\cth(s)\neq 0$, then $y'$ is bounded away from zero at least for $s$ big enough, hence $y$ blows up. Finally we have to study the special case of $\lim_{s\rightarrow +\infty}\cth(s)= 0$. We claim that, under this hypothesis, $z$ diverges. In fact by Corollary \ref{vth monotona} we have that $\lim_{s\rightarrow +\infty}\vth(s)$ is a zero of the function $f(\vth)=(\vthz-\lambda-\vth)\sth$. Therefore if $\cth$ converges to zero, necessarily $\lim_{s\rightarrow +\infty}\vth(s)=\vthz-\lambda$ and, by what we proved so far, $z$ diverges. By \eqref{eee} we get that
\begin{equation}\label{eq03}
y'=e^z\cth=\frac{\lambda\cth}{\vth-\vthz+\lambda}.
\end{equation}
Since $\cos(\vthz-\lambda)= 0$,we have 
$$
\lim_{s\rightarrow +\infty}\frac{\cth(s)}{\vth(s)-\vthz+\lambda}=\lim_{\vth\rightarrow\vthz-\lambda}\frac{\cth}{\vth-\vthz+\lambda}=\pm 1.
$$
Therefore by \eqref{eq03} $y'$ is bounded away from zero with a sign depending on $\vthz$ and $\lambda$, hence $y$ blows up in this case too. \cvd

\begin{figure}[H]
\centering
\includegraphics[width=0.9\textwidth]{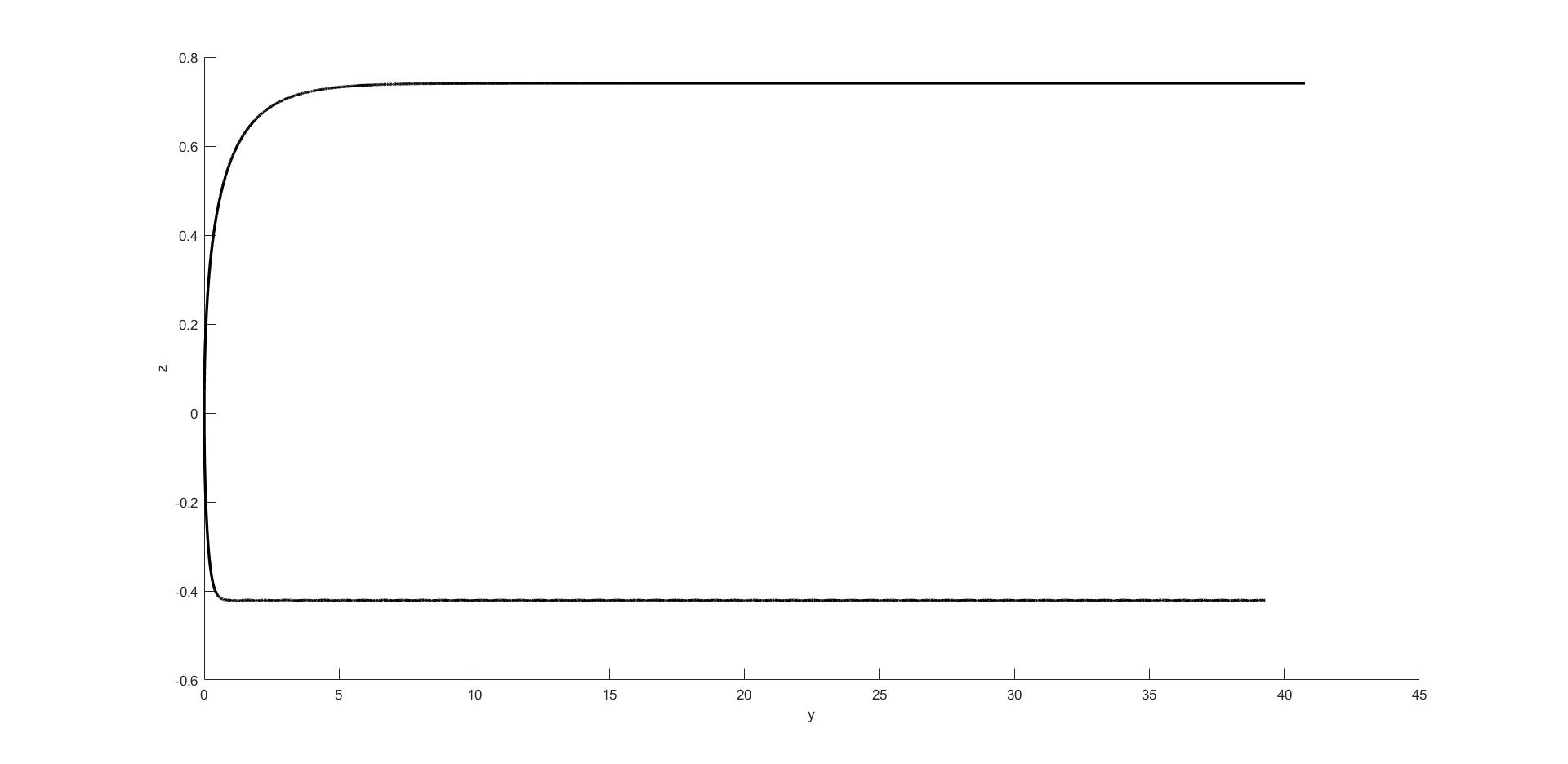}
\caption{The tilted grim reaper in $\sss$: $X=F_1$, $V=3F_2$ and $\vthz=\frac{\pi}{2}$.}
\label{fig01}
\end{figure}

\begin{figure}[H]
\centering
\includegraphics[width=0.9\textwidth]{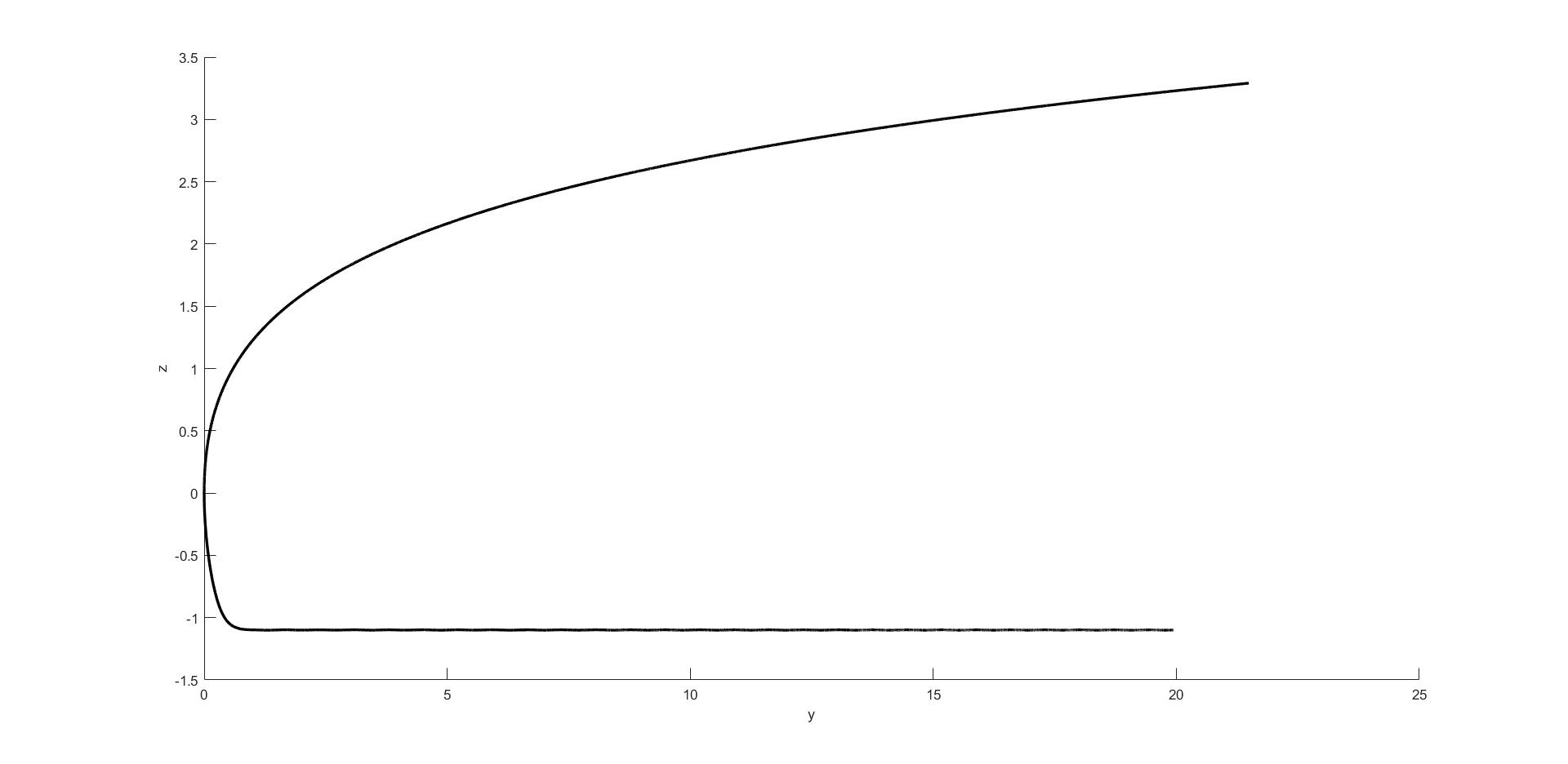}
\caption{Example of graphical translator defined on a half-plane with a critical point for $y$: $X=F_1$, $V=\frac{\pi}{4}F_2$ and $\vthz=\frac{\pi}{2}$.}
\label{fig02}
\end{figure}

\begin{figure}[H]
\centering
\includegraphics[width=0.9\textwidth]{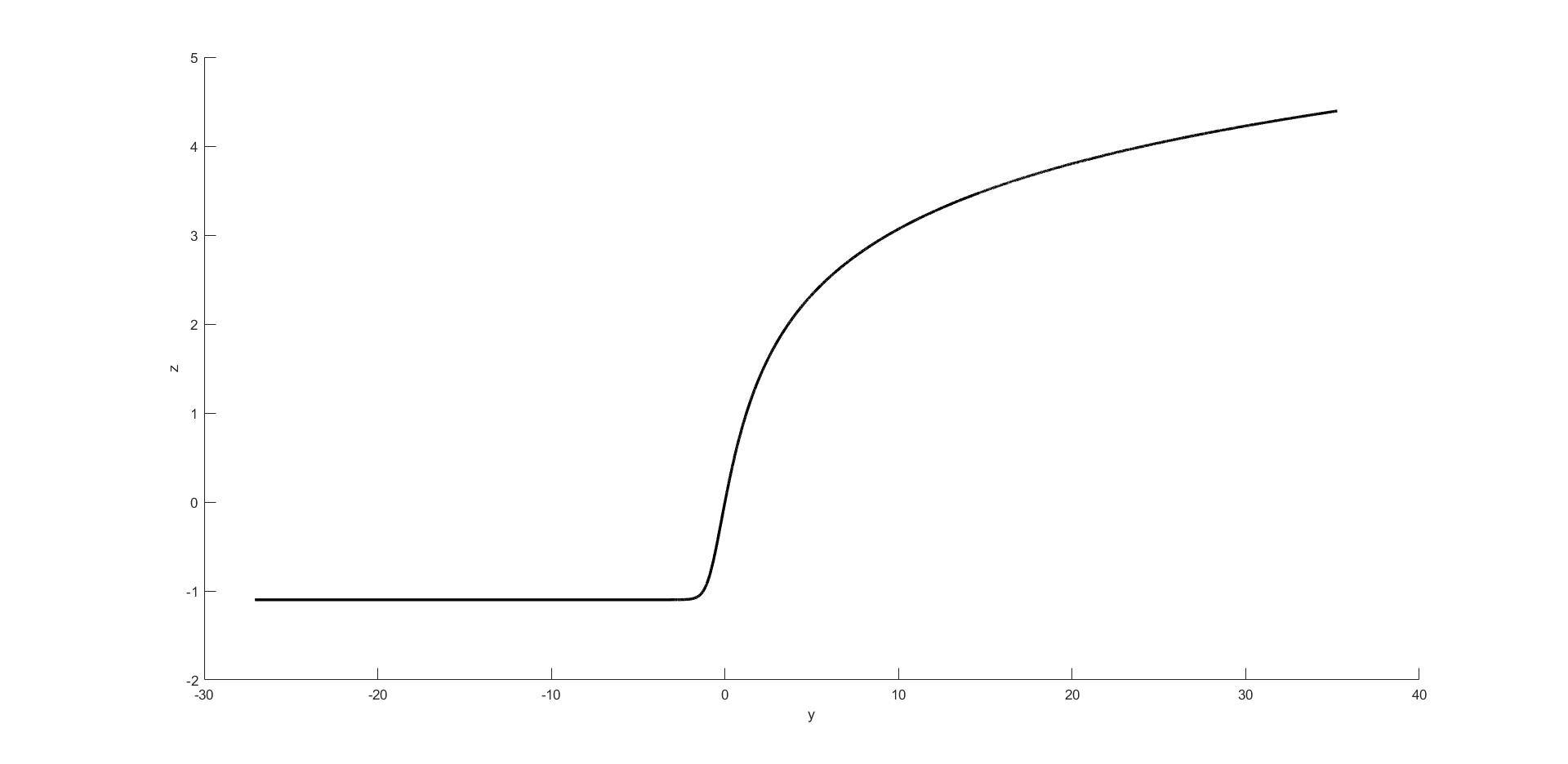}
\caption{Example of graphical translator defined on a half-plane without critical point for $y$: $X=F_1$, $V=-\frac{\pi}{8} F_2$ and $\vthz=\frac{\pi}{4}$.}
\label{fig03}
\end{figure}

\begin{figure}[H]
\centering
\includegraphics[width=0.9\textwidth]{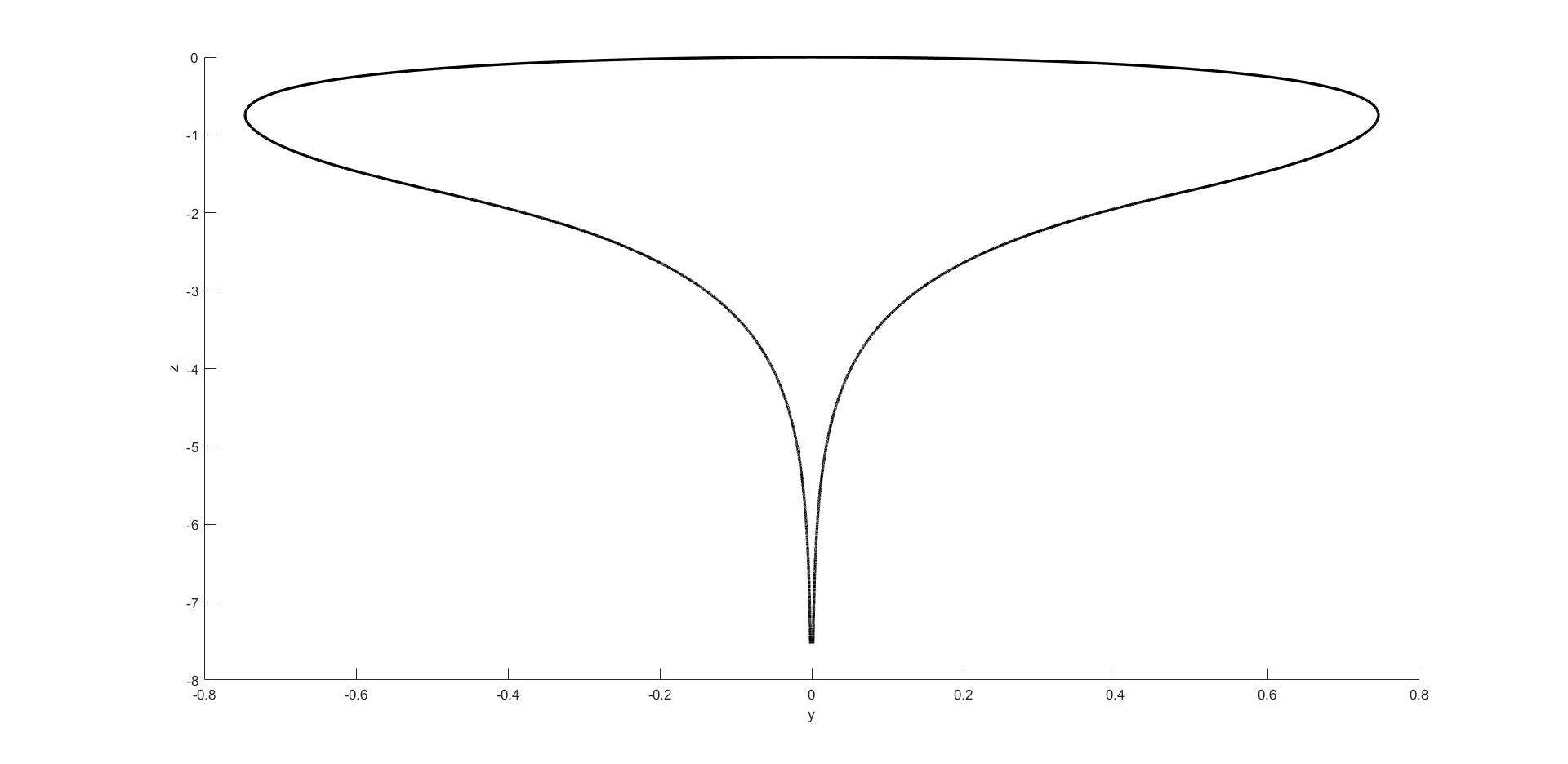}
\caption{Example with $z\rightarrow-\infty$ in both directions, symmetric case: $X=F_1$, $V=-F_3$ and $\vthz=0$. There are two critical point for $y$ and one for $z$.}
\label{fig04}
\end{figure}

\begin{figure}[H]
\centering
\includegraphics[width=0.9\textwidth]{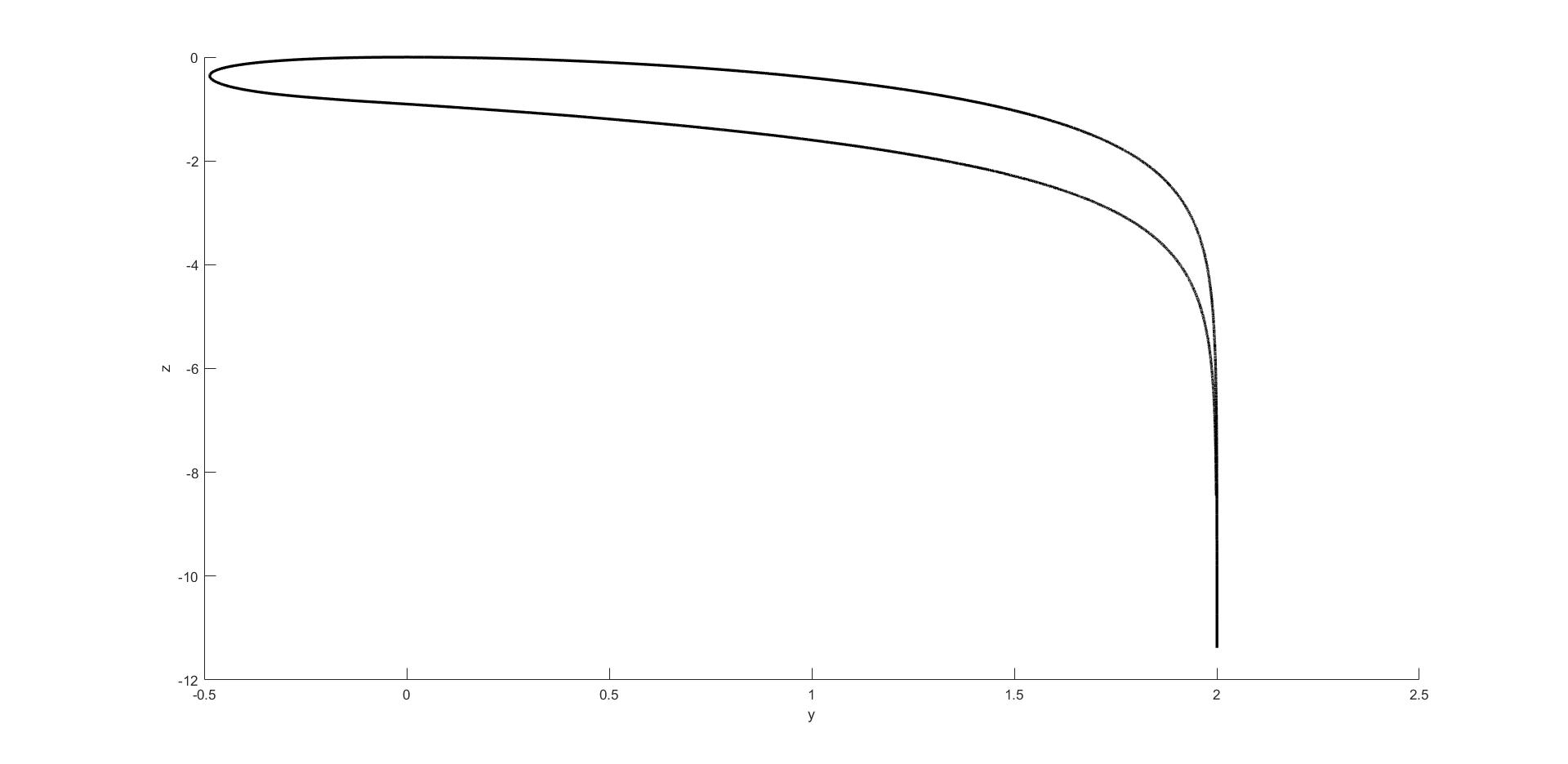}
\caption{Example with $z\rightarrow-\infty$ in both directions, non symmetric case: $X=F_1$, $V=2F_2-F_3$ and $\vthz=0$. There are two critical point for $y$ and one for $z$.}
\label{fig05}
\end{figure}

\begin{figure}[H]
\centering
\includegraphics[width=0.9\textwidth]{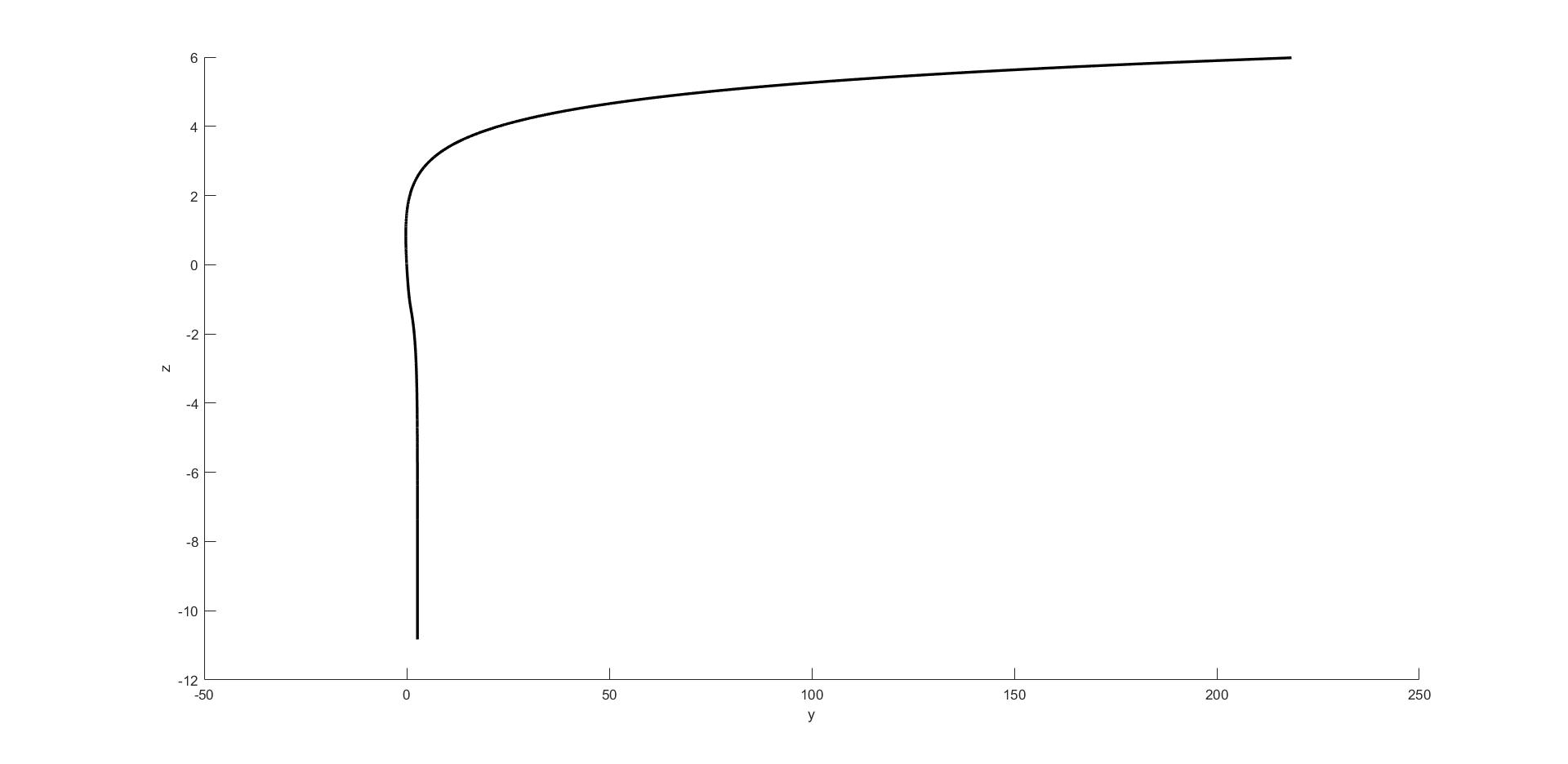}
\caption{Example with $z$ unbounded in both directions: $X=F_1$, $V=\frac 45 F_2-\frac{3}{10}F_3$ and $\vthz=2$. There is a critical point for $y$ and $z$ is strictly monotone.}
\label{fig06}
\end{figure}

\begin{figure}[H]
\centering
\includegraphics[width=0.9\textwidth]{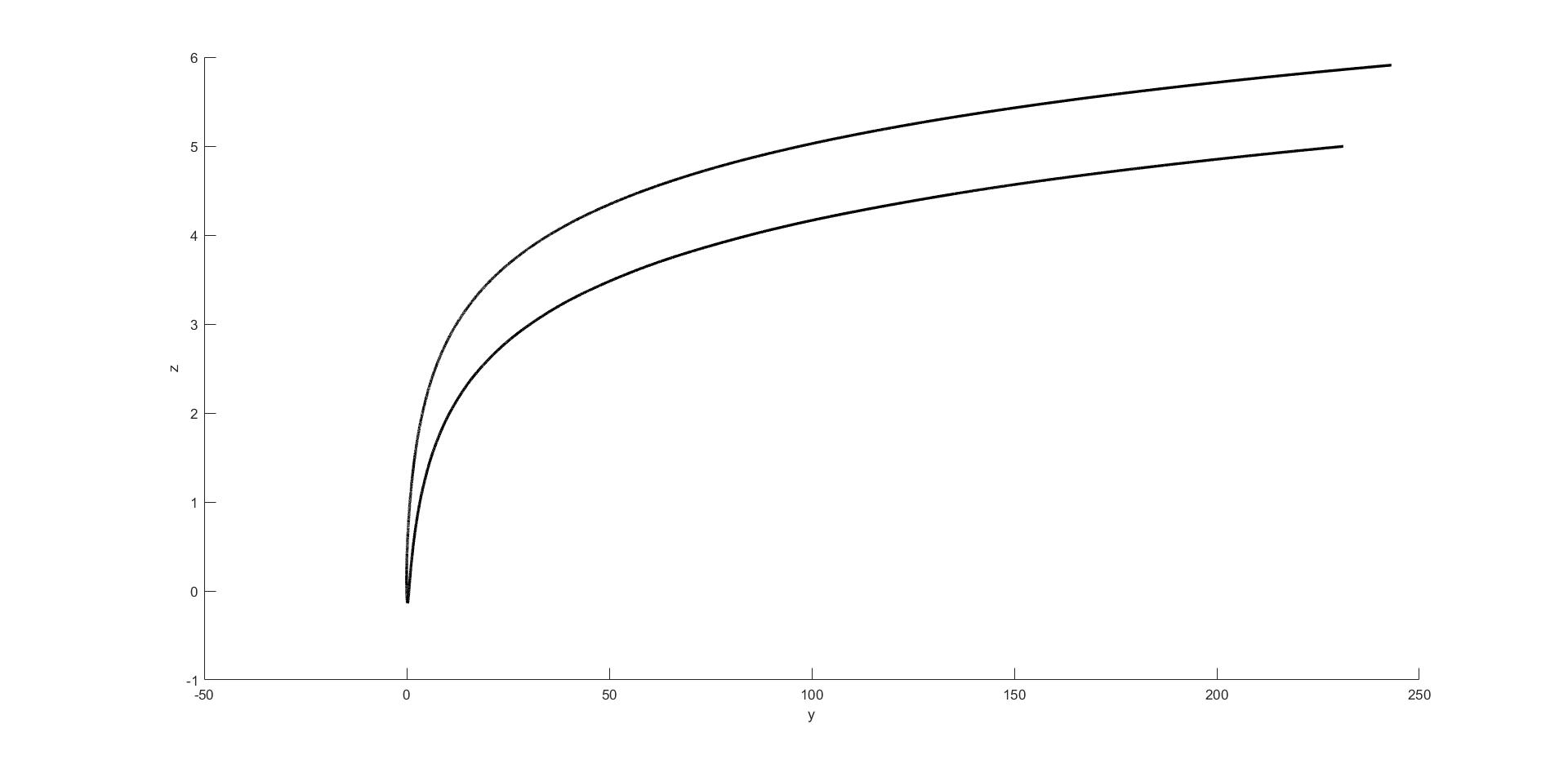}
\caption{Example with $z\rightarrow+\infty$ in both directions: $X=F_1$, $V=3F_2+3F_3$ and $\vthz=2$.}
\label{fig07}
\end{figure}

\begin{Remark}\label{rmk01}
\begin{itemize}
\item[1)] As a consequence of the previous result at least one of the two limits $\lim_{s\rightarrow\pm\infty}z(s)$ is finite, in fact at least one of the two limits $\lim_{s\rightarrow\pm\infty}\vth(s)\neq\vthz-\lambda$.
\item[2)] Since $z$ is strictly monotone in this case, $\mm$ can be globally parametrized as a graph $y=F(z)$ for some function $F$. It follows that in $\sss$ there are graphical translators that are defined on a half-plane, see Figures \ref{fig02} and \ref{fig03} for examples. This phenomenon cannot occur in the Euclidean space as showed by Shahriyari in \cite{Sh}.
\item[3)] If $z$ converges to a finite constant in both directions, then by Theorem \ref{main F2} there is a $k\in\mathbb Z$ such that $\vth$ varies in $(k\pi,(k+1)\pi)$. Therefore by \eqref{sistema01} $y$ has a critical point and it is unique by Lemma \ref{punti critici}. It follows that in this case $y$ has a global minimum (resp. maximum) and diverges to $+\infty$ (resp. $-\infty$) in both directions. These examples can be thought as the analogous in $\sss$ of the tilted grim-reaper cylinders in the Euclidean space. See figure \ref{fig01} for an example.
\item[4)] The Gaussian and the intrinsic sectional curvature of a $F_1$-invariant surface in $\sss$ have been computed in \cite{LM2}:
$$
\det(Ag^{-1})=-(\vth'+\cth)\cth,\quad K=-\sthq-\cth\vth'.
$$
By \eqref{vth4} we have that 
\begin{eqnarray*}
K & =  (\vth-\vthz+\lambda)\sth\cth-\sthq.
\end{eqnarray*}
Hence we can see that both the Gaussian curvature and $K$ can change their sign on $\mm$. The same phenomenon holds, for example, in the Heisenberg group \cite{Pi}. It is in contrast with the rigidity of the translators in the Euclidean space: in \cite{SX} the authors proved that any mean-convex and complete graphical translator in $\mathbb R^3$ is convex. 
\end{itemize}
\end{Remark}

We can to improve Lemma \ref{main F2} describing better the asymptotic shape of $\gamma$ and $\mm$. The curve $\gamma$ has two ends: let $\gamma^+=\gamma\cap\left\{z>0\right\}$ and let $\gamma^-=\gamma\cap\left\{z<0\right\}$. Analogously we define $\mm^+=\mm\cap\left\{z>0\right\}$ and let $\mm^-=\mm\cap\left\{z<0\right\}$. By Remark \ref{rmk01} we have that $\gamma^-$ is always asymptotic to a horizontal line, hence $\mm^-$ is asymptotic to a horizontal plane $z=z_0$ for some constant $z_0<0$. For the upper part we have different possibilities collected in the following result.

\begin{Lemma}\label{f2f2}
Fix $\mu=0$, $\lambda\neq 0$ and $\vthz\in\mathbb R$. 
\begin{itemize}
\item[1)] If $-\pi<\lambda<\pi$ and $\vthz=\lambda+k\pi$ for some $k\in\mathbb Z$, then $\mm^+$ is asymptotic to a half-logarithmic surface \ref{half}.
\item[2)] If $-\frac{\pi}{2}<\lambda<\frac{\pi}{2}$ and $\vthz=\lambda+\frac{\pi}{2}+k\pi$ for some $k\in\mathbb Z$, then $\mm^+$ is asymptotic to a plane $y=C_1z+C_0$ for some constant $C_0,C_1$.
\item[3)] In all other cases $\mm^+$ is asymptotic to a horizontal plane or a minimal logarithmic surface \eqref{log}.
\end{itemize}
\end{Lemma}
\proof We will focus on the curve $\gamma$, after that by $F_1$-invariance we can recover the analogous properties for $\mm$.
\begin{itemize}
\item[1)] Arguing by periodicity and changing the orientation of $\gamma$ if necessary, \eqref{vth4} can be reduced to 
$$
\left\{\begin{array}{rcl}
\vth' & = & -\vth\sth,\\
\vthz & \in & (-\pi,\pi)\backslash\left\{0\right\}.
\end{array}\right.
$$
Let us consider the case $\vthz>0$, the other one is analogous. In this case we have that $\lim_{s\rightarrow\infty}\vth(s)=0$. Therefore for $s$ big enough we have that
$$
\vth'\sim-\vth^2,\quad z'\sim\vth,\quad y'\sim e^{z}.
$$
Hence
$$
\vth\sim \frac 1s,\quad z\sim\log s,\quad y'\sim \frac{s^2}{2}.
$$
Combining the asymptotic expansion of $y$ and $z$ it is evident that $\mm^+$ is asymptotic to a half-logarithmic surface \eqref{half}.
\item[2)] By periodicity and changing the sign of $\vth$ if $k$ is odd, we can reduce \eqref{vth4} to
$$
\left\{\begin{array}{rcl}
\vth' & = & \left(\frac{\pi}{2}-\vth\right)\sth,\\
\vthz & \in & (0,\pi)\backslash\left\{\frac{\pi}{2}\right\}.
\end{array}\right.
$$
We get the end $\gamma^+$ when $\lim_{|s|\rightarrow+\infty}\vth(s)=\frac{\pi}{2}$. It follows that
$$
\vth'\sim\frac{\pi}{2}-\vth,\quad z'\sim 1.
$$
Therefore
$$
\vth\sim\frac{\pi}{2}+Ce^{-s},\quad z\sim s,
$$
for some constant $C$.
Hence
$$
y'=e^z\cth\sim -e^s\sin(Ce^{-s})\sim-C.
$$
Then both $y$ and $z$ grow linearly and $\mm^+$ is asymptotic to a plane.
\item[3)] Up to reverse the orientation of $\gamma$, we can suppose that $\gamma^+=\left\{\gamma(s)\ |\ s>0\right\}$. Suppose now that $z$ is unbounded, then by Lemma \ref{main F2} and \eqref{eee} we have that $$\lim_{s\rightarrow+\infty}z(s)=+\infty,\quad\lim_{s\rightarrow+\infty}\vth(s)=\vthz-\lambda.
$$
Let us define $\sigma_1:=\sin(\vthz-\lambda)$ and $\sigma_2:=\cos(\vthz-\lambda)$. Under our hypothesis $\sigma_1\in(0,1)$. Therefore for $s$ big enough we get
$$
z'\sim\sigma_1,\quad\Rightarrow\quad z\sim\sigma_1 s,
$$
and
$$
y'\sim\sigma_2 e^{\sigma_1s},\quad\Rightarrow\quad y\sim\frac{\sigma_2}{\sigma_1}e^{\sigma_1s}.
$$
The thesis follows easily. 
\end{itemize} \cvd

We conclude this section with some numerical simulations. 

\subsection{$F_1$-invariant translators in a general direction}
In this section we complete the proof of Theorem \ref{main1} describing the behavior of $\gamma$ for a general direction of evolution.

\begin{Lemma}
For any value of $\mu\neq 0$, $\lambda$ and $\vthz$ we have that, as $|s|$ diverges, the function $z$ diverges. In particular if $z\rightarrow -\infty$, then $y$ converges to the constant $-\frac{\lambda}{\mu}$, if $z\rightarrow +\infty$ $y$ can either converge to the constant $-\frac{\lambda}{\mu}$ or diverge. When both $y$ and $z$ diverge, $\gamma$ is asymptotic to a logarithmic minimal surface \eqref{log}.
\end{Lemma}
\proof
We consider the case of $s\rightarrow+\infty$, the other one is analogous. We know by Corollary \ref{vth monotona} that $\vth$ is converging to a zero of $f=\mu\cth+(\vthz-\lambda-\vth)\sth$. Since $\mu\neq 0$, then $\sth$ is converging to a non zero constant. By \eqref{sistema01} $z'=\sth$, we get that $z$ diverges linearly for any choice of $V$ with $\mu\neq 0$. We recall the formula \eqref{vth3} that links the angle with the coordinate funtions:
$$
\vth-\vthz+\lambda = e^{-z}(\lambda+\mu y).
$$
As in the proof of Theorem \ref{main F2}, we can deduce a lot of information from the fact that the left hand side of the previous equation is converging to a finite constant $k$. It is easily to see that, if $z$ is going to $-\infty$, then $y$ is converging to the constant $-\frac{\lambda}{\mu}$. On the other hand, when $z$ goes to $+\infty$ and $k\neq 0$, then $y$ necessarily blows up exponentially with a sign that depends on $k$ and $\mu$. Finally suppose that $z\rightarrow+\infty$ and $k=0$. The condition $k=0$ is satisfied if and only if $f(\vthz-\lambda)=0$, therefore we can deduce that $\cos(\vthz-\lambda)=0$. Moreover in this case $\lambda\neq 0$ otherwise $f(\vthz)=0$. By \eqref{vth3} we have
$$
y'=(\lambda+\mu y)\frac{\cth}{\vth-\vthz+\lambda}.
$$
Note that under our hypothesis $z'=\sth>0$ at least for $s$ big enough, then $\frac{\cth}{\vth-\vthz+\lambda}$ necessarily converges to $-1$. It follows that for $|s|$ big enough 
$$
y\sim-(\lambda+\mu y).
$$
Therefore in this case too $y$ can converge to $-\frac{\lambda}{\mu}$ or it can blow up exponentially with a sign that depends on $\lambda$ and $\mu$.
\cvd

\noindent We conclude this section showing some numerical simulations. 

\section{Translator invariant with respect to a generic direction}\label{SEZB}
Let us consider a generic Killing vector field $X=aF_1+bF_2+cF_3$. As usual we denote with $X$ also the one-parameter group of isometries generated by this vector field. We want to study the $X$-invariant translators: their properties are very different depending whether $c=0$ or not. We start considering $c=0$ proving Theorem \ref{main2}. We suppose that $ab\neq 0$,  otherwise we are in the situation described in Theorem \ref{main1}. Moreover it is not restrictive to assume $a=1$. Analogously to what we have seen in \eqref{rrr}, we can parametrize $\mm$ as
\begin{equation}\label{surfG}
\mm(u,s)=\left(u,bu+y(s),z(s)\right).
\end{equation}
As before we can introduce the angle function $\vth$, hence we have
\begin{equation}\label{sistG}
\left\{\begin{array}{rcl}
y'& =& e^z\cth ,\\
z' & = & \sth,
\end{array}\right.
\end{equation}
and, since $\sss$ is homogeneous, we can suppose that $y(0)=z(0)=0$. Since $b\neq 0$, the equation for $\vth$ is more complicated than in \cite{LM}, in particular we are not able to reproduce the strategy that we used to transform \eqref{vth2} into \eqref{vth4}. In order to compute $\vth'$, first we need to know the mean curvature of $\mm$.

\begin{Lemma}\label{mcG}
Let $\mm$ as in \eqref{surfG}, then its mean curvature is
$$
H=\frac{e^{2z}\left(\left(e^{4z}+b^2\right)\vth'-2b^2\cth\sthq\right)}{\left(e^{4z}+b^2\sthq\right)^{\frac 32}}.
$$
\end{Lemma}
\proof
From \eqref{surfG}, it is easy to see that a basis of tangent vector field is 
\begin{eqnarray*}
T_u:=\frac{\partial}{\partial u}\mm & = & F_1+bF_2 \ =\ e^zE_1+be^{-z}E_2,\\
T_s:=\frac{\partial}{\partial s}\mm & = & \cth E_2+\sth E_3.
\end{eqnarray*}
Therefore the induced metric on $\mm$ and a unit normal vector field are respectively
\begin{equation}\label{1ff}
g=\left(\begin{array}{cc}
e^{2z}+b^2e^{-2z} & be^{-z}\cth\\
be^{-z}\cth & 1
\end{array}\right),
\end{equation}
\begin{equation}\label{nu}
\nu=\frac{1}{\left(e^{4z}+b^2\sthq\right)^{\frac 12}}\left(b\sth E_1-e^{2z}\sth E_2+e^{2z}\cth E_3\right).
\end{equation}
By \eqref{LC} we have
\begin{eqnarray*}
\bar \nabla _{T_u}T_u & = & \left(b^2e^{-2z}-e^{2z}\right)E_3,\\
\bar \nabla _{T_s}T_s & = & \left(\vth'+\cth\right)\left(-\sth E_2+\cth E_3\right),\\
\bar \nabla _{T_u}T_s & = & \bar \nabla _{T_s}T_u \ = \ e^z\sth E_1-be^{-z}\sth E_2+be^{-z}\cth E_3.
\end{eqnarray*}
It follows that the second fundamental form of $\mm$ is
\begin{equation}\label{2ff}
A=\frac{1}{\left(e^{4z}+b^2\sthq\right)^{\frac 12}}\left(\begin{array}{cc}
\left(b^2-e^{4z}\right)\cth & b\left(1+\sthq\right)e^z\\
b\left(1+\sthq\right)e^z & \left(\vth'+\cth\right)e^{2z}
\end{array}\right).
\end{equation}
Since the mean curvature is the trace of $Ag^{-1}$, after some straightforward computations we get the thesis. \cvd

Now we have all the ingredients for describing the main properties of the functions $\vth$, $y$ and $z$. Therefore by \eqref{surfG} we recover the properties of $\mm$ proving Theorem \ref{main2}.

\begin{Lemma}
Let $\mm$ defined in \eqref{surfG}, and let $V$ be a Killing vector field. If $\mm$ is a translator in the direction of $V$, then we have two possibilities:
\begin{itemize}
\item[1)] the function $\vth$ is constant, in this case $\mm$ is the horizontal plane $z=0$, hence it is minimal and $V$ is tangent to $\mm$,
\item[2)] the function $\vth$ is non-constant, but $\sup\vth-\inf\vth\leq\pi$. In this case $\mm$ is complete and embedded, $z$ is strictly monotone and $y$ can have at most one critical point. It follows that both $y$ and $z$ have a limit (finite or not) as $|s|$ diverges.
\end{itemize}
\end{Lemma}
\proof
Let as usual $V=\eta F_1+\lambda F_2+\mu F_3$ be a Killing vector field. Without loss of generality we can assume $\eta=0$, in fact if it is not the case we have that
$$
V=\eta T_u+\tilde{\lambda}F_2+\mu F_3,
$$
where $\tilde{\lambda}=\lambda-b\eta$. In \eqref{nu} we computed the expression of the normal $\nu$, therefore we get that $\mm$ is a translator in the direction of $V$ if and only if
\begin{eqnarray}
\nonumber H & = & \bg (\nu,V)\ = \ \bg(\nu,{\lambda}F_2+\mu F_3)\\
 & = &  \frac{e^z}{\left(e^{4z}+b^2\sthq\right)^{\frac 12}}\left(-{\lambda}\sth+\mu\left(e^z\cth-y\sth\right)-2\mu bu\sth \right)
\label{HG}\end{eqnarray}
Because of the symmetry of $\mm$, the last term in the equalities above cannot depend on $u$, it follows that either $\mu=0$ or $\sth \equiv 0$. In the latter case, by \eqref{sistG}, we have that $z\equiv 0$ and it is well known that this surface is minimal in $\sss$. We conclude this proof by considering the case $\mu=0$. By Lemma \ref{mcG} and \eqref{HG}, after some algebraic manipulation we have that
\begin{equation}\label{vthG}
\vth' = \frac{\sth}{e^{4z}+b^2}\left[b^2\left(2\sth\cth-{\lambda}e^{-z}\sthq\right)-{\lambda}e^{3z}\right].
\end{equation}
Unfortunately we are unable to release this equation from the others, as was done in Section \ref{SEZA}, but \eqref{vthG} is enough to deduce some information about $\mm$. Let $\vth_0=\vth(0)$ and suppose that $\sin\vth_0\neq 0$. First of all we can notice that $\vth$ is bounded: in fact for every $k\in\mathbb Z$ the constant function $\vth(s)=k\pi$ is a solution of \eqref{vthG} and it acts as a barrier, then for every $\vth_0$, there is a $k\in\mathbb Z$ such that for every $s$
\begin{equation}\label{GGG}
k\pi<\vth(s)<(k+1)\pi.
\end{equation}
In particular we have that $\vth$ is defined for every $s\in\mathbb R$. By \eqref{sistG}, the same holds for $y$ and $z$. Hence $\mm$ is complete. Finally by \eqref{GGG} and \eqref{sistG}, we can see that $z$ has no critical points, hence $\mm$ is embedded, and $y$ has at most one critical point. \cvd

Finally we prove Theorem \ref{main3}.

\noindent\emph{Proof of Theorem \ref{main3}.} Let $X=aF_1+bF_2+cF_3$ with $c\neq0$. Without loss of generalization we can suppose $c=1$. Let $\mm$ be a $X$-invariant surface, therefore by Lemma \ref{sottogruppi} there is a curve $\gamma(s)=(x(s),y(s),0)$ parametrized by arc-length such that $\mm$ can be parametrized as follows
$$
\begin{array}{rcl}
\mm(u,s) & = & \left(a(1-e^{-u}),b(e^{u}-1),u\right)\star\gamma(s)\\
& = & \left(a(1-e^{-u})+e^{-u}x(s),b(e^{u}-1)+e^{u}y(s),u\right).
\end{array}
$$
A basis of tangent vector field is 
\begin{eqnarray*}
T_u:=\frac{\partial}{\partial u}\mm(u,s) & = & aF_1+bF_2+F_3\\
& = &(a-x(s))E_1+(b+y(s))E_2+E_3,\\ 
T_s:=\frac{\partial}{\partial s}\mm(u,s) & = &x'(s)E_1+y'(s)E_2.
\end{eqnarray*}
It follows that a normal vector field is
$$
\nu=-y'(s)E_1+x'(s)E_2-\left[y'(s)(x(s)-a)+x'(s)(y(s)+b)\right]E_3.
$$
Let $V=\eta F_1+\lambda F_2+\mu F_3$ be a Killing vector field. Note that we can always suppose that $\mu=0$ because, since $c\neq0$ we can write
$$
V=\frac{\mu}{c}T_u+\tilde{\eta}F_1+\tilde{\lambda}F_2,
$$
where $\tilde{\eta}=\eta-a$ and $\tilde{\lambda}=\lambda-b$. Therefore $\mm$ is a translator in the direction of $V$ if and only if
\begin{equation*}
H=\bar g\left(\frac{\nu}{|\nu|},V\right)=\bar g\left(\frac{\nu}{|\nu|},\tilde{\eta}F_1+\tilde{\lambda}F_2\right).
\end{equation*}
After some straightforward calculation we have
\begin{equation}\label{ultima}
|\nu|H=-\tilde{\eta}y'(s)e^{u}+\tilde{\lambda}x'(s)e^{-u}.
\end{equation}
Because of the symmetry of $\mm$, the left hand term in the equation above depends only on $s$ and not on $u$, therefore the only possibility is that the right hand term vanishes, that is $H=0$ and $V$ is tangent to $\mm$. Finally, if $\mm$ satisfies hypothesis $2)$, then it is trivial that it is a translator in the direction of $V$.
\cvd

\begin{Remark}\label{rmkfin}
By Theorem \ref{main3} and equation \eqref{ultima}, we have that when $c\neq 0$, $\mm$ is a $X$-invariant translator in the direction of $V$ if and only if
$$
\tilde{\lambda}x'=\tilde{\eta}y'=0.
$$
Clearly $x'$ and $y'$ cannot be both zero. Therefore we deduce that such $\mm$ exists if and only if $\tilde{\lambda}\tilde{\eta}=0$.
Moreover if in this case $\tilde{\lambda}\neq 0$ (resp. $\tilde{\eta}\neq 0$) then $x'\equiv 0$ (resp. $y'\equiv 0$). It the special case  $a=b=0$ this means that the only $F_3$-invariant translators are the planes $x\equiv c_0$ and $y\equiv c_0$ for some constant $c_0$.
\end{Remark}

\bigskip

\noindent Giuseppe Pipoli, \emph{Department of Information Engineering, Computer Science and Mathematics, Universit\`a degli Studi dell'Aquila}, via Vetoio 1, 67100, L'Aquila, Italy.\\
Email: giuseppe.pipoli@univaq.it
\end{document}